\documentclass[11pt]{article}
\usepackage{amsmath}
\usepackage{amssymb}
\usepackage{graphicx}

\newcommand{\qed}{\hfill$\Box$\par\medskip\par\relax}

\numberwithin{equation}{section}

\newcommand{\eps}{\varepsilon}
\newcommand{\Z}{{\mathbb Z}}
\newcommand{\N}{{\mathbb N}}

\newcommand{\G}{{\mathcal G}}
\newcommand{\C}{{\mathcal C}}

\renewcommand{\b}{{\mathfrak b}}
\renewcommand{\phi}{\varphi}

\newcommand{\Cov}{{\mathop{\rm Cov}}}

\newcommand{\EV}{\mathbf{E}}
\newcommand{\PV}{\mathbf{P}}
\newcommand{\IE}{\mathbb{E}}
\newcommand{\IP}{\mathbb{P}}
\newcommand{\hIE}{\widehat{\mathbb{E}}}
\newcommand{\hIP}{\widehat{\mathbb{P}}}

\newcommand{\Pomega}{{\mathtt P}_{\!\ome}}
\newcommand{\Po}{{\mathtt P}_{\!\omega}}

\newcommand{\Eomega}{{\mathtt E}_{\ome}}

\newcommand{\ii}{{\hat\iota}}
\newcommand{\EEE}{\mathcal{E}}

\newcommand{\dg}{\mathop{\mathrm{deg}}}
\newcommand{\0}{\mathbf{o}}
\newcommand{\mm}{\mathfrak{m}}

\newcommand{\x}{\text{\boldmath ${\xi}$}}
\newcommand{\ome}{\text{\boldmath ${\omega}$}}
\newcommand{\T}{\text{\boldmath $T$}}
\newcommand{\bi}[1]{\stackrel{\leftrightarrow}{#1}}
\newcommand{\lp}[1]{\stackrel{\leftarrow}{#1}}
\newcommand{\rp}[1]{\stackrel{\rightarrow}{#1}}
\newcommand{\1}[1]{{\mathbf{1}}{\{#1\}}}

\allowdisplaybreaks

\newtheorem{theo}{Theorem}[section]
\newtheorem{lmm}[theo]{Lemma}

\newtheorem{prop}[theo]{Proposition}

\newtheorem{rmk}{Remark}[section]

\newcommand{\supp}{\mathop{\rm supp}}

\title{Random walks on Galton--Watson trees with random 
conductances}

\author{Nina Gantert$^{1}$ \and Sebastian M\"{u}ller$^{2}$ \and
 Serguei~Popov$^{3}$ \and Marina Vachkovskaia$^{3}$}

\begin{document}

\maketitle

{\footnotesize \noindent $^{~1}$
Technische Universit\"at M\"unchen, Fakult\"at f\"ur Mathematik,
Boltzmannstr. 3, 85748 Garching,
Germany\\
\noindent e-mail: \texttt{gantert@ma.tum.de}\\
\noindent 
url:\texttt{http://www-m14.ma.tum.de/en/staff/gantert/}

\noindent $^{~2}$ LATP, CMI Universit\'e de Provence
39 rue Joliot Curie, 13453 Marseille cedex 13, France\\
\noindent e-mail: \texttt{mueller@cmi.univ-mrs.fr},\\ 
\noindent 
url:\texttt{http://www.latp.univ-mrs.fr/$\sim$mueller/}

\noindent $^{~3}$Department of Statistics, 
Institute of Mathematics, Statistics and Scientific Computation,\\
University of Campinas--UNICAMP, 
rua S\'ergio Buarque de Holanda 651, 13083--859, Campinas SP,
Brazil\\
\noindent e-mail: \texttt{\{popov,marinav\}@ime.unicamp.br}, 
\newline \noindent
url: \texttt{http://www.ime.unicamp.br/$\sim$\{popov,marinav\}}
}

\begin{abstract}
We consider the random conductance model, where the underlying 
graph is an infinite supercritical Galton--Watson tree, 
the conductances are independent  but their distribution may depend on 
the degree of the incident vertices. 
We prove that, if the mean conductance is finite,
 there is a deterministic, strictly positive speed $v$ 
such that~$\lim_{n\to\infty} \frac{|X_n|}{n}= v$ a.s.\
(here,~$|\cdot|$ stands for the distance from the root).
We give a formula for~$v$ in terms of the laws of certain effective
conductances and show that, if the conductances share the same
expected value, the speed is not larger than the speed of simple
random walk on  Galton--Watson trees. The proof relies
on finding a reversible measure for the environment observed by the
particle.
\\[.3cm]\textbf{Keywords:} rate of escape,
 environment observed by the particle,
        effective conductance, reversibility
\\[.3cm]\textbf{AMS 2000 subject classifications:}
 60K37, 60J10
\end{abstract}

\section{Introduction}
This paper is a contribution to the theory of
random walks on random networks.
Here, the underlying graph is an infinite supercritical Galton--Watson
tree with independent conductances whose distribution may depend
on the degree of the incident vertices. It is not difficult to see
that such random walks are transient; see
Proposition~\ref{p_transience}. We denote the random walk by~$\{X_n\}_{n\in\N}$.  We say that there is a law
of large numbers if there exists a deterministic~$v$ (the rate
of escape, or the speed) 
such that~$\lim_{n\to\infty} \frac{|X_n|}{n}=v$ a.s.,\
where,~$|\cdot|$ stands for the distance from the root.  A standard method to prove  laws of large numbers  is to work in
the space of rooted weighted trees and to consider the environment
observed by the particle. 
This approach has the advantage, provided one is able to construct a
stationary measure, that it gives rise to a stationary ergodic
Markov chain and one can apply the ergodic theorem. We identify
the reversible measure for the environment in Section~\ref{s_env}
and prove a formula for the speed which involves effective
conductances of subtrees, 
see Theorem~\ref{thm:speed}. A first consequence is that the speed is
a.s.\ positive. For the case of non-degenerate random
conductances having the same mean
 we show a \emph{slowdown result}: the speed of the
random walk with random conductances is strictly smaller than the
speed of the simple random walk. 
Finally, we consider an example on the binary tree, see Proposition~\ref{ex1},  where explicit asymptotic
results are obtained. This example illustrates how the choice of
the random environment influences the speed of the random walk.

Simple random walks on Galton--Watson trees were studied
in~\cite{LPP:95} where among other results a law of large number is
proved, using the environment observed by the particle.
In~\cite{LP:09} one finds more references and details about
this and related models. There are mainly two generalizations of
this model. The first is the so-called~$\lambda$-biased random
walk. In this model the random walk chooses the direction towards
the root with probability proportional to~$\lambda$ while the
probability to choose any of the sites in the opposite direction
is proportional to~$1$. In~\cite{L:90} it was proved that
the~$\lambda$-biased random walk is positive recurrent
if~$\lambda> m$, null recurrent if~$\lambda= m$, and transient
otherwise. Here,~$m$ is the mean number of offspring of the
Galton--Watson process. In the transient case, it was shown
in~\cite{LPP:95} and~\cite{LPP:96} that~$|X_n|/n\to v_\lambda>0$ 
a.s.,
where~$v_\lambda$ is deterministic. An explicit formula 
for~$v_\lambda$ is
only known for~$\lambda=1$ (that is, for the case of SRW).
For~$\lambda\leq m$,~\cite{PZ:08} proves  a quenched central limit
theorem for~$|X_n| - n v$  by constructing a stationary
measure for the environment process. In the
critical case, $\lambda= m$, the  central limit theorem has the following form:
for almost every realization of the tree, the ratio $|X_{[nt]}|/\sqrt{n}$
converges in law as~$n\to\infty$ to a deterministic multiple
of the absolute value of a Brownian motion.
The second generalization are random walks in random environment
(RWRE) on Galton--Watson trees. The main difference to our work
is that while in our model the conductances are 
realizations of an independent environment, 
in the RWRE model the \textit{ratios} of
the conductances are realizations of an i.i.d.\ environment.
Therefore,
the behaviour of RWRE is richer; the walk may be recurrent or
transient and the speed positive or zero. We refer
to~\cite{Aid:08} and~\cite{Faraud} and references therein
for recent results. 

Our model can also be seen from a more general point of view as an example
of a stationary random network. A stationary random network is a random
rooted network whose distribution is invariant under re-rooting along the
path of the random walk (defined through the corresponding electric network)
started at the original root. This notion generalizes the concept of
transitive networks where the condition of transitivity  is replaced by the
assumption that an invariant distribution along the path of the random walk
exists. Under first moment conditions, this model is also known as
a unimodular random network, see \cite{AL:07} and \cite{BC}, or an invariant 
measure of a graphed equivalence
relation. In fact, unimodular random
networks correspond to stationary and
reversible random networks. A straightforward consequence of the
stationarity and the sub-additive ergodic theory, see e.g. \cite{AL:07}, is
the existence of the speed,  i.e., for almost every realization of a
stationary and reversible random network, $ |X_n|/n$ converges almost surely.

The rest of the paper is organized as follows. In
Section~\ref{sec_model} we give a formal description and notations
of the model. The environment observed by the
particle is introduced in Section~\ref{s_env} and in
Section~\ref{s_speed} we present the main results that are proved
in Section~\ref{s_proofs}. Some
open questions are in Section~\ref{s_open}.

\section{The model}\label{sec_model}
A rooted tree~$\T$ is a nonoriented, connected, and  locally
finite graph without loops. One vertex~$\0$ is singled out and
called the root of the tree. The rooted tree is then denoted
by~$(\T,\0)$. We use the same notation~$\T$ for the set of
vertices of the tree and the tree itself; the set of edges is
denoted by~$\EEE(\T)$. For a vertex~$x\in\T$  we denote
by~$\dg(x)$ the degree of~$x$ (i.e., the number of edges incident
to~$x$). The \emph{index} of~$x$ is defined by $\ii(x)=\dg(x)-1$.
Let~$|x|$ be the (graph) distance from~$x$ to the root. We
write~$x\sim y$ if ~$x$ and~$y$ are connected by an edge,
i.e.,~$(x,y)\in\EEE(\T)$. 
Then, for a fixed tree~$\T$ and any nonnegative integers~$k,m$,
define
\begin{align*}
  U_{k,m}(\T) = & \{(x,y) \in\EEE(\T) :  \ii(x)=k, \ii(y)=m\}
\end{align*}
to be the set of edges connecting vertices of indices~$k$ and~$m$.
An electrical network is a graph where each edge has a positive
label called the \emph{conductance} or \emph{weight} of the edge.
In our model these conductances are realizations of 
a collection of independent
random variables. More precisely, for every unordered
pair~$\{k,m\}$ we label all edges~$e\in U_{k,m}$ with positive
i.i.d.\ random variables $\xi(e)$ with common
law~${\tilde\mu}_{k,m}$.  We
denote by~$\gamma_{k,m}$ the expected value of~$\xi$
under~${\tilde\mu}_{k,m}$ (note that $\gamma_{k,m}\in (0,\infty]$ 
for all~$k,m$)
and write~$\x:=(\xi(e),~e\in\EEE(\T))$ for the
environment of conductances (weights) on the tree. Clearly, the
above definitions are symmetric in the sense
that~$U_{k,m}(\T)=U_{m,k}(\T)$,
${\tilde\mu}_{k,m}={\tilde\mu}_{m,k}$,
$\gamma_{k,m}=\gamma_{m,k}$ for all $k,m$. Such a weighted
rooted tree is then denoted by the triple~$(\T,\0,\x)$.

Now, we would like to consider a model where the tree itself is
chosen at random. Let~$p_0, p_1, p_2, p_3, \ldots$ be the
parameters of a Galton--Watson branching process, i.e., $p_k$ is
the probability that a vertex has~$k$ descendants. We assume that
$p_0=0$, see Remark
\ref{rem_cond_survival} for the case where this assumption is
dropped.
 Furthermore, suppose that
\begin{equation}
 \label{sum_pj}
\mu: = \sum_{j=1}^\infty jp_j\in (1, +\infty).
\end{equation}
The latter  guarantees that there exists~$j>1$ such that~$p_j>0$,
so that the tree a.s.\ has infinitely many ends.

\begin{figure}
\centering
\includegraphics{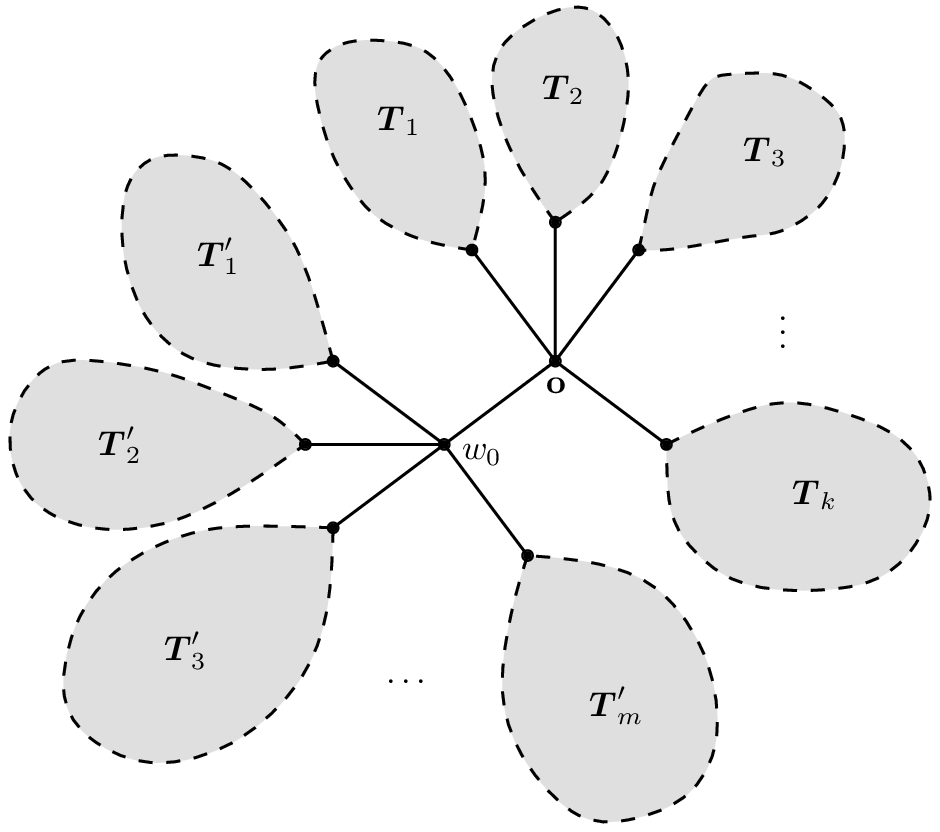}
\caption{On the definition of~$\IP_{k,m}$: $\T_1,\ldots,\T_k$ and
$\T'_1,\ldots,\T'_m$ are i.i.d.\ weighted Galton--Watson  trees
with law~$\hIP$} \label{f_Ekm}
\end{figure}

\begin{figure}
\centering
\includegraphics{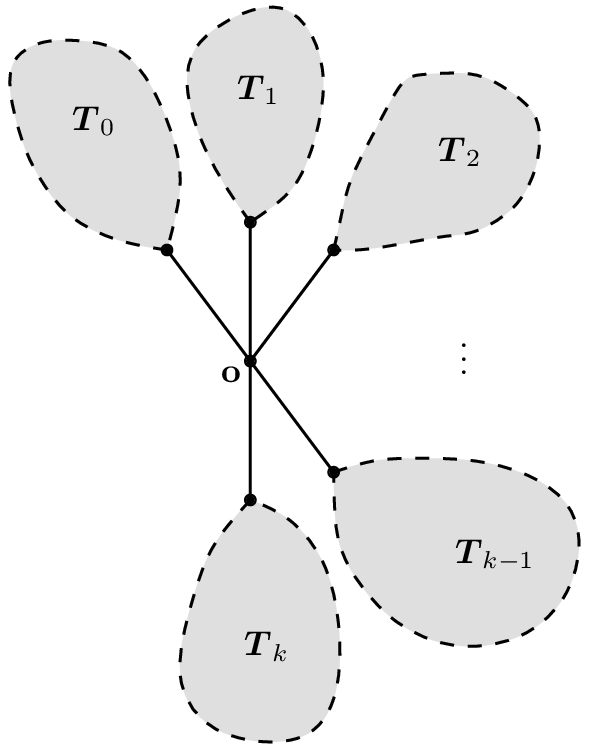}
\caption{On the definition of~$\IP_k$:~$\T_0,\ldots,\T_k$ are
i.i.d.\ weighted Galton--Watson trees with law~$\hIP$}
\label{f_Ek}
\end{figure}

Define~$\hIP,~\hIE$ to be the probability and expectation for the
usual rooted Galton--Watson tree (i.e., the genealogical tree of
the Galton--Watson process with the above parameters) with random
conductances as  described above. 
Now define $\IP_{k,m}$ in the following way, see Figure~\ref{f_Ekm}.
Take i.i.d.\ copies $\T_1,\ldots,\T_k,\T'_1,\ldots,\T'_m$ 
 of a weighted Galton--Watson tree with law~$\hIP$. Denote the roots
of
 $\T_1,\ldots,\T_k$ by $w_1,\ldots,w_k$.
Take a vertex $\0$ with $\ii(\0) = k$ and attach vertices
$w_1,\ldots,w_k$ with edges  $\ell_1,\ldots,\ell_k$, 
starting from~$\0$.
In the same way, attach $\T'_1,\ldots,\T'_m$ to edges starting 
from a second vertex 
$w_0$. Choose the conductances of all this edges 
independently according to the corresponding laws. Finally, 
connect~$o$ and $w_0$ by an edge and choose its conductance
independently from everything according to ${\tilde\mu}_{k,m}$.
 We denote by $\IE_{k,m}$ the expectation with respect to
$\IP_{k,m}$.
For each $k$, we can now define 
$\IP_{k} = \sum_{m=1}^\infty p_m \IP_{k,m}$, and we denote by
$\IE_k$ its expectation.
Note that $\IP_{k}$ is the law of the weighted
Galton--Watson tree,  conditioned on the
event~$\{\ii(\0)=k\}$, see Figure~\ref{f_Ek}.
Note that with this construction, under~$\IP_k$ the subtrees 
attached to $w_0,\ldots,w_k$ are independent
and have the law~$\hIP$.

The probability measure~$\IP$ for the
\emph{augmented} Galton--Watson tree with conductances is given by the mixture $\IP=\sum_{k=1}^\infty p_k\IP_k$. In other words, first we choose an index~$k$ with probability~$p_k$, and
then sample the random tree from the measure~$\IP_k$. We note
that this is equivalent to considering two independent weighted
Galton--Watson trees with law~$\hIP$ connected by a weighted
edge;  the conductance of this edge is sampled from the
corresponding distribution independently of everything. 
We denote the corresponding expectation by~$\IE$.
The
important advantage of considering augmented weighted
Galton--Watson trees is the following stationarity property: for
any non-negative functions~$f,g,u$ on the space of rooted weighted trees we have
\begin{align}
\lefteqn{ \IE \big[ f(\T,\0,\x)g(\T,w_0(\0),\x)
u(\xi(\ell_0(\0)))\big]}\nonumber\\
& = \IE \big[g(\T,\0,\x)f(\T,w_0(\0),\x) u(\xi(\ell_0(\0)))\big].
\label{AGW_stationary}
\end{align}
Indeed, using the representation of~$\IE_{k,m}$ shown in
Figure~\ref{f_Ekm}, it is straightforward to obtain that
\begin{align*}
 \lefteqn{\IE \big[ f(\T,\0,\x)g(\T,w_0(\0),\x) u(\xi(\ell_0(\0)))
\big]}\\
 &= \sum_{k,m} p_k p_m \IE_{k,m} \big[ f(\T,\0,\x)g(\T,w_0(\0),\x) 
u(\xi(\ell_0(\0)))\big]\\
 &= \sum_{m,k} p_m p_k \IE_{m,k} \big[g(\T,\0,\x)f(\T,w_0(\0),\x)) 
u(\xi(\ell_0(\0))\big]\\
 &= \IE \big[g(\T,\0,\x)f(\T,w_0(\0),\x) u(\xi(\ell_0(\0)))\big].
\end{align*}

Let us denote
\begin{equation}\label{rever}
\pi_x=\sum_{z\sim x}\xi(x,z)
\end{equation}
and  define the discrete time random walk~$\{X_n\}_{n\in \N}$
on~$\T$ in the environment~$\ome=(\T,\0,\x)$ through the
transition probabilities
\[
 q_\ome(x,y) = \frac{\xi(x,y)}{\pi_x}.
\]
For a fixed realization~$\ome$ of the environment,  denote
by~$\Pomega, \Eomega$ the probability and expectation with respect
to the random walk~$\{X_n\}_{n\in \N}$, so that
$q_\ome(x,y)=\Pomega[X_{n+1}=y\mid X_n=x]$ and $\Pomega[X_0=\0]=1$.
The definition~(\ref{rever}) 
implies that this random walk is reversible
with the corresponding reversible measure~$\pi$, that
is, for all~$x,y\in\T$ we
have~$\pi_xq_\ome(x,y)=\pi_yq_\ome(y,x)=\xi(x,y)$.

It is not difficult to obtain that the random walk defined
above is a.s.\ transient:
\begin{prop}
\label{p_transience}
 The random walk $\{X_n\}_{n\in \N}$ is transient for $\IP$-almost 
all environments~$\ome$.
\end{prop}

\noindent \textit{Proof.} The random walk is transient if and only 
if the effective conductance of the tree (from the root to infinity)
is strictly positive, see Theorem~2.3 of~\cite{LP:09}.
By~\eqref{sum_pj}, we can
choose~$\delta$ and~$d$ such that
\begin{equation}\label{deltaso}
(1-\delta)\sum_{i=1}^d jp_j>1.
\end{equation}
Then, choose~$\eps$ small enough such that
\[
 {\tilde\mu}_{k,m}[(\eps, \infty)]\ge 1-\delta
\]
for all~$k,m \le d$. We define a percolation process on $\T$ by
deleting all edges with~$\xi(e)\le \eps$. This process dominates a
Bernoulli percolation on a~$d+1$-regular tree with retention
parameter~$1-\delta$. Due to \eqref{deltaso} this
percolation process is supercritical. Hence, there is a.s.\ 
an infinite subtree of the original tree (not necessarily containing
the root) such that all the
conductances of this subtree are at least~$\eps$. 
Since this subtree is itself an infinite Galton--Watson tree, the random walk on
it is transient and it has positive effective conductance. We
conclude that also the effective conductance of the original tree is
strictly positive.
 \qed

\begin{rmk}
Under the condition that $\gamma = \sum_{k,m} p_k p_m \gamma_{k,m}<\infty$, Proposition \ref{p_transience} is a special case of Proposition 4.10 in \cite{AL:07}. 
\end{rmk}

\section{Environment observed by the particle}
\label{s_env} The aim of this section is to construct a reversible
measure for the environment, observed by the particle.

Let~$\gamma = \sum_{k,m} p_k p_m \gamma_{k,m}$, and define
\begin{equation}
\label{def_m(T)}
 \mm(\T,\0,\x) = \frac{\pi_{\0}}{\ii(\0)+1}
\end{equation}
Loosely speaking,~$\mm$ is the mean conductance from $\0$ to its
neighbours. Clearly, we have
\[
 \IE\left(\mm(\T,\0,\x)\right) = \sum_k
\frac{p_k}{k+1}\IE_k(\pi_{\0})
  = \sum_{k,j} p_k p_j \gamma_{k,j} = \gamma.
\]
Provided that $\gamma < \infty$, we can define a new probability
measure~$\PV$ on the set of
weighted rooted trees through the corresponding expectation
\begin{equation}
\label{def_E_v}
 \EV\big[ f(\T,\0,\x)\big] = \frac{1}{\gamma} \IE\big[\mm(\T,\0,\x)
f(\T,\0,\x)\big].
\end{equation}
Also, for two~$\PV$-square-integrable functions~$f,g$, we define
their scalar product
\begin{equation}
\label{def_scalar}
 (f,g) = \EV \big[f(\T,\0,\x)g(\T,\0,\x)\big].
\end{equation}

The environment observed by the particle is the
process on the space of all weighted rooted trees with transition
operator
\begin{align}
 G f(\T,\0,\x) &= \sum_{z\sim\0} q_\ome(\0,z) f(\T,z,\x)\nonumber\\
  &= \frac{1}{\pi_{\0}}  \sum_{z\sim\0} \xi(\0,z) f(\T,z,\x).
\label{def_G}
\end{align}
Let us now prove that~$G$ is reversible with respect to~$\PV$. In
particular, this implies that~$\PV$ is a stationary measure for
the environment, observed by the particle.
\begin{lmm}
\label{l_revers} For any two functions~$f,g\in L_2(\PV)$, we have
$(f,Gg)=(Gf,g)$.
\end{lmm}

\noindent \textit{Proof.} Indeed, we have
\begin{align}
 (f,Gg) &= \frac{1}{\gamma}
  \IE\Big[\frac{1}{\ii(\0)+1}f(\T,\0,\x)
           \sum_{z\sim\0} \xi(\0,z) g(\T,z,\x)\Big]\nonumber\\
 &=\frac{1}{\gamma}\sum_k \frac{p_k}{k+1}
  \IE_k\Big[f(\T,\0,\x)\sum_{j=0}^k\xi(\ell_j(\0)) g(\T,w_j(\0),\x)
\Big]
 \nonumber\\
 &=\frac{1}{\gamma}\sum_k p_k
 \IE_k\big[f(\T,\0,\x)\xi(\ell_0(\0)) g(\T,w_0(\0),\x)\big]
\nonumber\\
 &= \frac{1}{\gamma}\IE\big[f(\T,\0,\x)\xi(\ell_0(\0))
 g(\T,w_0(\0),\x)\big]. \label{conta_revers1}
\end{align}
In the same way we obtain
\begin{equation}
\label{conta_revers2}
 (g,Gf) = \frac{1}{\gamma}\IE\big[g(\T,\0,\x)\xi(\ell_0(\0))
 f(\T,w_0(\0),\x)\big],
\end{equation}
and so, using~\eqref{AGW_stationary},
we conclude the proof of Lemma~\ref{l_revers}. \qed

\section{Main results}
\label{s_speed} 
Usually, for any weighted rooted tree~$(\T,\0,\x)$ we will write
just~$\T$ since it is always clear from the context to which root
and to which set of weights we are referring. Let~$\C(\T)$ be the
effective conductance from the root to infinity (cf.\ e.g.\
Section~2.2 of~\cite{LP:09}). Suppose that the random walk starts
at the root, i.e.,~$X_0=\0$. Provided that the following limit exists, we
define the speed of the random
walk~$\{X_n\}_{n\in \N}$ by
\begin{equation}
\label{def_speed}
 v =\lim_{n\to\infty}\frac{|X_n|}{n}\, .
\end{equation}

\begin{figure}[h]
\centering
\includegraphics{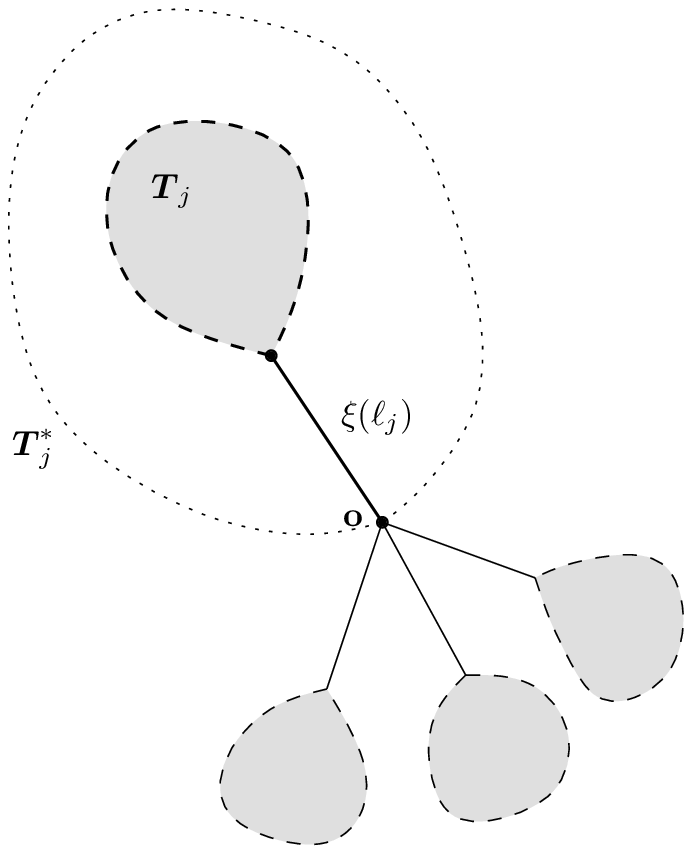}
\caption{Definition of the tree~$\T_j^*$}~\label{f_T_star}
\end{figure}

Recall that the neighbours of the root~$\0$ are denoted by~$w_0,
\ldots, w_{\ii(\0)}$, while $\ell_0, \ldots,
\ell_{\ii(\0)}$ are the corresponding edges. Denote
$\xi_j:=\xi(\ell_j)$. Let~$\T_j$ be the subtree of~$\T$ rooted 
at~$w_j$ and~$\T_j^*$ be the tree~$\T_j$ together with the 
edge~$\ell_j$ (see Figure~\ref{f_T_star}; we assume that the 
root of~$\T_j^*$ is~$\0$). Note also that
\begin{equation}
\label{calc_cond}
 \C(\T_j^*)=\frac{1}{\frac{1}{\xi_j}+\frac{1}{\C(\T_j)}}.
\end{equation}

One of the main results of this paper is the following formula for
the speed
of the random walk with random conductances:
\begin{theo}
\label{thm:speed} 
Assume $\gamma < \infty$. Then, the limit in~\eqref{def_speed}
exists~$\Pomega$-a.s. for $\IP$-almost all $\omega$. Moreover, $v$ is
deterministic and is given by
\begin{align}
 v &=1-  \frac2\gamma \IE\Big(\xi_0\frac{\C(\T_0^*)}{\C(\T)}\Big)
\label{eq:speed'}\\
 &= \sum_{k=1}^\infty p_k\Big[1-\frac{2}{\gamma}
  \IE_k
\Big(\xi_0\frac{\C(\T_0^*)}{\C(\T)}\Big)\Big]\label{eq:speed}\\
&=\sum_{k=1}^\infty p_k\Big[1-\frac{2}{(k+1)\gamma}
  \IE_k \Big(\sum_{i=0}^k \xi_i\frac{\C(\T_i^*)}{\C(\T)}\Big)\Big].
\label{eq:speed''}
\end{align}
\end{theo}

\begin{rmk}
\label{rem_cond_survival}
We can also consider the case when $p_0>0$, i.e., when the
augmented Galton--Watson process may die out. In this case, we
have to condition on the survival of the process. We then
 obtain the following formula:
\begin{align}
\lim_{n\to \infty} \frac{|X_n|}n &=1-  \frac2\gamma
\IE\Big(\xi_0\frac{\C(\T_0^*)}{\C(\T)}\Big|survival\Big)\\
&=\sum_{k=1}^\infty p_k
\frac{1-q^{k+1}}{1-q^2}\Big[1-\frac{2}{\gamma}
  \IE_k
\Big(\xi_0\frac{\C(\T_0^*)}{\C(\T)}\Big)\Big],
\end{align} ~$\Pomega$-a.s. for $\IP$-almost all $\omega$, where $q$
is the extinction probability of the
Galton--Watson process. The relation of the latter formulas with 
(\ref{eq:speed'}) and (\ref{eq:speed''}) is the same as
in~\cite{LPP:95} for simple random walk.
\end{rmk}

From~\eqref{eq:speed'}--\eqref{eq:speed''} it is not immediately
clear if the speed is positive, so let us prove the following
\begin{theo}
\label{pos_speed}
Assume that $\gamma < \infty$. Then, the quantity~$v$ given
in~\eqref{eq:speed'} is strictly positive.
\end{theo}

\begin{rmk}
In the case of bounded conductances, i.e., if there exists $c,C>0$
such that
$\supp {\tilde\mu}_{k,m} \subseteq [c, C]$,
Theorem~\ref{pos_speed} also follows from~\cite{Virag:00} and the
fact that
supercritical Galton--Watson trees have the anchored expansion
property, see~\cite{CP:04}.
\end{rmk}

Next, we treat also the case where the expected conductance in some
edges may be infinite:
\begin{theo}
\label{zero_speed}
Assume that there exist $k, m$ such that $\gamma_{k,m}= \infty$.
Then, the limit in~\eqref{def_speed} is $0$,~$\Pomega$-a.s. for
$\IP$-almost all $\omega$.
\end{theo}

Using Theorem~\ref{thm:speed}, we can compare the speed of the
random walk on Galton--Watson trees with random conductances to the
speed of simple random walk (SRW) on the same tree (observe that SRW
corresponds to the
case when all the conductances are a.s.\ equal to the same
positive constant). Let~$v_{SRW}$ be the speed of SRW on the
Galton--Watson tree; by Theorem~3.2 of~\cite{LPP:95} it holds that
\begin{equation}
\label{speed_SRW}
 v_{SRW} = \sum_k p_k\frac{k-1}{k+1}.
\end{equation}

\begin{theo}
\label{compara_speeds}
Assume $\gamma < \infty$. 
Let~$v$ be the speed of the random walk $\{X_n\}_{n\in \N}$.
\begin{itemize}
\item[(i)] We have
\begin{equation}
\label{speed_cov}
 v = v_{SRW} - \frac{2}{\gamma}
\Cov\Big(\xi_0,\frac{\C(\T^*_0)}{\C(\T)}\Big)
\end{equation}
(the covariance is with respect to~$\IE$). 
\item[(ii)] Suppose
that the conductances have the same expectation, 
i.e.,\ $\gamma_{k,m} = \gamma$, for all $k, m$, and~$\xi_0$ is a
non-degenerate random variable. Then
\begin{equation}
\label{smallerthansimple}
v< v_{SRW} .
\end{equation}
\end{itemize}
\end{theo}

In practice, it is not easy to use Theorem~\ref{thm:speed} for the
exact calculation of the speed due to the following reason. While
it is not difficult to write a distributional equation that the
law of~$\C(\T_0^*)$ should satisfy, it is in general not possible
to solve this equation explicitly. Nevertheless,
Theorem~\ref{thm:speed} can be useful, as the following example
shows. Let us consider the binary tree (i.e.,~$p_2=1$) with
i.i.d.\ conductances
\[
\xi=\begin{cases} 1,& \text{ with probability } 1-\eps_n, \\
 a_n, &   \text{ with probability } \eps_n,
           \end{cases}
\]
where~$\eps_n\to 0$ and~$a_n \to \infty$ as~$n\to \infty$.
Let~$v_n$ be the speed of the random walk with conductances
distributed as above.
\begin{prop}
\label{ex1} Assume that $\eps_n a_n\to\eta \in [0,\infty]$ as $n \to
\infty$.
Then,
\[
 \lim_{n\to\infty} v_n =\frac{1}{3(\eta+1)}=  \frac{v_{SRW}}{\eta+1}.
\]
\end{prop}

\section{Proofs}
\label{s_proofs}

\noindent \textit{Proof of Theorem~\ref{thm:speed}.} The first
part of the proof is to show ergodicity of our process. Since we
follow here
the arguments in~\cite{LPP:95}, see also~\cite{LP:09}
(Section~16.3), we only give a sketch. To make
use of the ergodic theorem it is convenient to  work on the space
of bi-infinite paths. A bi-infinite path
$\ldots,x_{-1},x_0,x_1,\ldots$ is denoted by~$\bi{x}$. 
We denote by~$\rp{x}$ the path~$x_0,x_1,\ldots$ and by~$\lp{x}$ the
path~$\ldots,x_{-1},x_0.$ The path of the random walk has the
property that it converges a.s.\ to a boundary point; this follows
from transience. The space of convergent paths~$\bi{x}$ in~$\T$ is
denoted by~$\bi{\T}$ (convergent means here that one has convergence
both for $n \to \infty$ and $ n \to -\infty$). We consider the
(bi-infinite) path space
\[
 \textsf{PathsInTrees}:=\left\{(\bi{x},\T): \bi{x}\in\bi{\T} 
          \right\}.
\]
The rooted tree corresponding to~$(\bi{x},\T)$ is~$(\T,x_0)$.
Define the shift map:
\[
 (S\!\bi{x})_n:=x_{n+1}, \quad
S(\bi{x},\T):=(S\! \bi{x},\T)
\] 
and write $S^k$ for the $k$th iteration.
In order to define a probability measure on \textsf{PathsInTrees}
we extend the random walk to all integers by letting~$\lp{x}$ be
an independent copy of~$\rp{x}$. We use the notation~$RW\times
\PV$ for the corresponding measure on \textsf{PathsInTrees}.
 Observe that due to the reversibility of the probability measure $\PV$
(see Lemma~\ref{l_revers}), the corresponding Markov chain,
describing the environment and the path seen from the current
position of the walker,  
is stationary. We proceed by a regeneration argument. Define the set
of regeneration points
\[
 \textsf{Regen}:=\{(\bi{x},\T)\in
\textsf{PathsInTrees}: x_{-n}\neq
x_0 \text{ and } 
x_n\neq x_{-1}\text{ for all } n>0\}.
\]

The first step is to show that a.s.\ the trajectory has infinitely many
regeneration points. To this end, define the set of ``fresh''
points:
\[
 \textsf{Fresh}:=\{ (\bi{x},\T)\in \textsf{PathsInTrees}:
x_n\neq x_0\text{ for all } n<0\}.
\]
The idea is to show that the trajectory of the particle a.s.\ has
 infinitely many fresh points, 
 and then one concludes by observing that a positive
fraction of the fresh points has a (uniform) positive probability
to be a regeneration point. The first fact follows from a.s.\
transience of the random walk and the fact that two independent random walks converge a.s.\ to different ends. Moreover, there exists a positive density of fresh points. To see
this, observe that the sequence $\big\{ \1
{S^n(\bi{x},\T)\in\textsf{Fresh}},~n\in\N \big\}$ is stationary and
hence $\frac1n \sum_{i=1}^n \1 {S^i(\bi{x},\T)\in\textsf{Fresh}}$
converges to some positive (random) number $b$. We define the (a.s.\
positive) random variable
\[
 U(\bi x,\T)=\min_{z\sim x_0} \Pomega
[X_1=z,~X_n\neq x_0 \text{ for all } n\geq 1].
\]
Again, the sequence $\{ \1 {U(S^n(\bi{x},\T))>\eps},~n\in\N  \}$ is
stationary for any $\eps>0$  and we have that  $\frac1n \sum_{i=1}^n \1{U(S^i(\bi{x},\T))>\eps}$ converges to some (random) $c(\eps)>0$.
For every realization of the process (in the bi-infinite path
space) we can choose $\eps$ sufficiently small in such a way that
$c(\eps)>1-b/2$. Therefore,
\[
 \liminf_{n\to\infty}\frac1n \sum_{i=1}^n
\1{S^i(\bi{x},\T)\in\textsf{Fresh}} \1
{U(S^i(\bi{x},\T))>\eps}>b/2.
\] 
Eventually, this shows the existence of a random sequence 
$(n_k)_{k\in\N}$
such that~$S^{n_k}(\bi{x}, \T)$ is a fresh point and
$U(S^{n_k}(\bi{x},\T))>\eps$ and hence there exist infinitely many
regeneration points. Again we  follow the arguments
in Section~16.3 of~\cite{LP:09}. Let $x$ be some vertex in $\T$. We denote by $\T_{x}$ the subtree of $\T$ formed by those edges that become disconnected from $\0$ when $x$ is removed.  Define $n_{\textsf{Regen}}=\inf\{n>0: S^{n}(\bi{x},\T)\in \textsf{Regen}\}$. To each $(\bi{x},\T)\in\textsf{Regen}$ we associate a so-called slab:
\[
\textsf{Slab}(\bi{x},\T)=(\langle x_{0}, x_{1},\ldots,x_{n-1}\rangle, \T\setminus (\T^{x_{-1}}\cup\T^{x_{n}})),
\] 
where $n=n_{\textsf{Regen}}$ and $\langle x_{0},x_{1},\ldots,x_{n-1}\rangle$  stands for the path of the walk from time $0$ to time $n-1$. Write $S_{\textsf{Regen}}=S^{n_{\textsf{Regen}}}$ when $(\bi{x},\T)\in \textsf{Regen}$ and consider the random variables $\textsf{Slab}(S^{k}_{\textsf{Regen}}(\bi{x},\T))$. In contrast to~\cite{LP:09} these random variables are not independent. However, if we define $\textsf{Index}(\bi{x},\T)=\ii(x_{0})$ then due to the construction of our model we have that $\textsf{Slab}(S^{k}_{\textsf{Regen}}(\bi{x},\T))$ conditioned on $\textsf{Index}(S^{k}_{\textsf{Regen}}(\bi{x},\T))$ is an independent sequence.
In order to obtain an i.i.d.\ sequence we denote by $\textbf{i}$ the smallest possible index, i.e., $\textbf{i}=\inf\{i\geq 1: p_{i}>0\}$ and define 
\[
 \textsf{Regen}_{\textbf{i}}:=\{(\bi{x},\T)\in
\textsf{Regen}: \ii(x_0)=\textbf{i}\}.
\]
Since $\textsf{Index}(S^{k}(\bi{x},\T))$ is a stationary Markov chain on $\{i:~p_{i}>0\}$ which is irreducible and recurrent one shows that  $\textsf{Index}(S^{k}_{\textsf{Regen}}(\bi{x},\T))$ is a recurrent Markov chain. To see this let us first treat the case where  $\{i:~p_{i}>0\}$ is finite. Then, the probability that $\textsf{Index}(S^{k+1}_{\textsf{Regen}}(\bi{x},\T))=\textbf{i}$ conditioned on $\textsf{Index}(S^{k}_{\textsf{Regen}}(\bi{x},\T))$ is a random variable bounded away from zero. For the general case, we proceed similarly to the proof that there is an infinite number of regeneration points. In fact, we show first that there is a positive fraction of regeneration times. Then, define $V(\bi{x},\T)=\ii(\0)$
and consider the stationary sequence $\{\1{V(S^{n}(\bi{x},\T)\}\leq K}, n\in\N\}$ for some $K\in\N$. Choose $K$ sufficiently large such that there are infinitely many regeneration times whose index is smaller than $K$ and proceed as in the finite case.

Eventually,  there is an infinite number of index $\textbf{i}$ regeneration points. Since $\textsf{Slab}(S^{k}_{\textsf{Regen}_{\textbf{i}}}(\bi{x},\T))$ is an i.i.d.\ sequence that generates the whole tree and the random walk, we obtain that the
system~$(\textsf{PathsInTrees}, RW\times \PV, S)$ is ergodic.

As in the proof of Theorem~3.2 in~\cite{LPP:95} we calculate the
speed as the increase of the \emph{horodistance} from a boundary
point. So let~$\b$ be a boundary point,~$x$ be a vertex in~$\T$,
and let us denote by~$\mathcal{R}(x,\b)$ the ray from~$x$ to~$\b$.
Given two
vertices we can define the confluent~$x\wedge_\b y$ with respect
to $\b$ as the vertex 
where the two rays~$\mathcal{R}(x,\b)$
and~$\mathcal{R}(y,\b)$ coalesce.
 We define the
signed distance from~$x$ to~$y$ as~$[y-x]_\b:=|y-x\wedge_\b y | -
|x-x\wedge_\b y |$. (Imagine, you sit in $\b$ and wonder how
many steps more you have to do to reach~$y$ than to reach~$x$.)
Denote by~$x_{-\infty}$ (respectively,~$x_{+\infty}$) the boundary
points towards which~$\lp{x}$ (respectively,~$\rp{x}$) converges.
Since~$x_{-\infty}\neq x_{+\infty}$ a.s., there exists some
constant~$c$ such that for all sufficiently
large~$n$ we have~$|x_n-x_0|=[x_n-x_0]_{x_{-\infty}} +c$. 
(More precisely, $c = 2 |x_0 - x_0 \wedge_{x_\infty} x_\infty|$.)
Now, the speed is the limit
\[
\lim_{n\to\infty} \frac1n [x_n-x_0]_{x_{-\infty}}=
\lim_{n\to\infty} \frac1n \sum_{k=0}^{n-1}
[x_{k+1}-x_k]_{x_{-\infty}}.
\]
Since~$(\textsf{PathsInTrees}, RW\times \PV, S)$ is ergodic, these
are averages over an ergodic stationary sequence, and hence by the
ergodic theorem converge a.s.\ to their mean
\begin{equation}
\label{eq_speed_ergodic}
 v =\int [x_1-x_0]_{x_{-\infty}} d (RW \times \PV)(\bi{x},\T).
\end{equation}
The remaining part of the proof is devoted to find a more explicit
expression for this mean. This step is more delicate in the present
situation than for
SRW. Recall that~$\lp{x}$ is an independent copy of~$\rp{x}$.
Hence, we are interested in the probability that a random walk steps
\emph{towards} the boundary point of a second independent random
walk. 

We say that the random walk~$\{X_n\}_{n\in \N}$ escapes to
infinity in the direction~$\ell_k$, if
\[
 |\T_k\cap \{X_0, X_1, X_2, \ldots\}|=\infty.
\]
Observe that transience implies that the random walk escapes to
infinity in only one direction (since otherwise~$\0$ would be 
visited infinitely many times). Let us define a random variable
$\Theta$ in the following way:~$\Theta=k$ iff the random walk
escapes to infinity in the direction~$\ell_k$. Let
\[
 \psi_\xi=\Pomega^{\0}[X'_1=w_{\Theta}(\0)]
\]
stand for the probability that an independent
copy~$\{X'_n\}_{n\in\N}$
of the random walk~$\{X_n\}_{n\in\N}$ makes the first step in
the escape direction of~$\{X_n\}_{n\in\N}$. Hence, we can write
equation~(\ref{eq_speed_ergodic}) as
\begin{equation}
\label{eq_V}
 v=-\EV \psi_\xi+ \EV (1-\psi_\xi)=1-2 \EV \psi_\xi.
\end{equation}

Let us compute~$\psi_\xi$ now.

\medskip
\noindent
\textbf{Claim.} We have
\begin{equation}
\label{klaro}
\Pomega[\Theta=k]=  \frac{\C(\T^*_k)}{\C(\T)}\, .
\end{equation}

\noindent
\textit{Proof of the claim.} This is, of course, a standard fact, but
we still write its proof for completeness. 
Let
\[
 \tau_y=\inf\{n: X_n=y\}.
\]
Denote by
\[
 \eta_x (y)=\Pomega^x[\tau_y=\infty]
\]
the probability that the random walk starting from~$x$ never
hits~$y$. Note that
\begin{equation}
 \label{esc1}
\eta_{w_k(\0)}(\0)=\frac{\C(\T_k^*)}{\xi_k}
 =\frac{1}{\xi_k}\cdot\frac{1}{\frac{1}{\xi_k}+\frac{1}{\C(\T_k)}}
 =\frac{\C(\T_k)}{\xi_k+\C(\T_k)}.
\end{equation}
This follows e.g.\
from formula~(2.4) of~\cite{LP:09} and the fact that,
 for the random walk with random
conductances restricted to~$\T_k^*$, the escape probability from
the root equals~$\eta_{w_k(\0)}(\0)$.

Due to the Markov property,
\begin{align*}
 \Pomega[\Theta=k]&
  =\frac{\xi_k}{\pi_{\0}}\big(\eta_{w_k(\0)}(\0)
       +(1-\eta_{w_k(\0)}(\0))\Pomega[\Theta=k]\big)\\
& \qquad {} + \sum_{j\ne
k}\frac{\xi_j}{\pi_{\0}}(1-\eta_{w_j(\0)}(\0))\Pomega[\Theta=k],\\
\end{align*}
so, using~\eqref{esc1}, we obtain
\begin{align*}
  \Pomega[\Theta=k]&=\Big(1-\sum_{j=0}^{\ii(\0)}
    \frac{\xi_j}{\pi_{\0}}\cdot
 \frac{\xi_j}{\xi_j+\C(\T_j)}\Big)^{-1}\frac{\xi_k}{\pi_{\0}}\cdot
 \frac{\C(\T_k)}{\xi_k+\C(\T_k)}\\
&=\frac{\frac{\xi_k}{\pi_{\0}}\cdot
    \frac{\C(\T_k)}{\xi_k+\C(\T_k)}}{\sum_{j=0}^{\ii(\0)}
    \frac{\xi_j}{\pi_{\0}}\cdot
 \frac{\C(\T_j)}{\xi_j+\C(\T_j)}}\\
& = \frac{\C(\T^*_k)}{\C(\T)},
\end{align*}
which finishes the proof of the claim.
\qed

\medskip

 Now, we have
\begin{align*}
 \psi_\xi=\sum_{k=0}^{\ii(\0)}\frac{\xi_k}{\pi_{\0}}
          \Pomega[\Theta=k].
\end{align*}
Then, using~\eqref{def_E_v} and plugging in~\eqref{klaro}, we have
\begin{align}
 \EV \psi_\xi&=\EV\Big(\pi_{\0}^{-1}
\frac{\sum_{k=0}^{\ii(\0)}\xi_k
\C(\T^*_k)}{\C(\T)}\Big)\nonumber\\
 &=\sum_{j=1}^\infty
\frac{p_j}{(j+1)\gamma}\IE_j\Big(\frac{\sum_{k=0}^j\xi_k
    \C(\T_k^*)}{\C(\T)}\Big)\nonumber\\
 &= \sum_{j=1}^\infty
\frac{p_j}{\gamma}\IE_j\Big(\xi_0\frac{\C(\T_0^*)}{\C(\T)}
\Big)\nonumber\\
&=\frac{1}{\gamma} \IE\Big(\xi_0\frac{\C(\T_0^*)}{\C(\T)}
\Big). \nonumber
\end{align}

Together with~\eqref{eq_V}, this
implies~\eqref{eq:speed'}, \eqref{eq:speed}, and~\eqref{eq:speed''}.
\qed

\noindent \textit{Proof of Theorem~\ref{pos_speed}}. For each~$j$,
let $Z_1^{(j)},Z_2^{(j)},Z_3^{(j)},\ldots$ be i.i.d.\ random
variables, having the distribution of the effective conductance of
the tree $\T_0^*$, conditioned on the event that the root has
index~$j$. Denote $r_{k,j}=p_kp_j\gamma_{k,j}/\gamma$; 
observe that
$\sum_{k,j}r_{k,j}=1$, and $r_{k,j}=r_{j,k}$.
Assume that $(Z_i^{(j)})_{i=1,2, \ldots }$ are independent
collections of random variables for $j=1,2, \ldots $.
Then we have
\begin{align*}
 &\sum_{k,j} r_{k,j} E\Big(\frac{Z_1^{(j)}+\cdots +Z_j^{(j)}}
  {Z_1^{(j)}+\cdots +Z_j^{(j)}
 +Z_1^{(k)}+\cdots +Z_k^{(k)}}\Big)\\
&~~=\sum_{k,j} r_{k,j} \Big(1-E\Big(\frac{Z_1^{(k)}+\cdots
+Z_k^{(k)}}{Z_1^{(j)}+\cdots +Z_j^{(j)}
 +Z_1^{(k)}+\cdots +Z_k^{(k)}}\Big)\Big)\\
&~~=1-\sum_{k,j} r_{k,j} E\Big(\frac{Z_1^{(k)}+\cdots
+Z_k^{(k)}}{Z_1^{(j)}+\cdots +Z_j^{(j)}
 +Z_1^{(k)}+\cdots +Z_k^{(k)}}\Big),
\end{align*}
so, by symmetry,
\begin{equation}
 \label{sum_pkpj}
\sum_{k,j} r_{k,j} E\Big(\frac{Z_1^{(j)}+\cdots
+Z_j^{(j)}}{Z_1^{(j)}+\cdots +Z_j^{(j)}
 +Z_1^{(k)}+\cdots +Z_k^{(k)}}\Big)=\frac{1}{2}.
\end{equation}

We have (one may find it helpful to
look at Figure~\ref{f_Ekm} again)
\begin{align*}
 \IE\Big(\xi_0 \frac{\C(\T_0^*)}{\C(\T)} \Big)
 &=\sum_{k,j} p_k p_j\IE_{k,j}\Big(\xi_0 \frac{\C(\T_0^*)}{\C(\T)}
\Big)\\
&=\sum_{k,j} p_k p_j \int\limits_0^\infty x\IE_{k,j} \Big(
\frac{\C(\T_0^*)}{\C(\T)} 
\mid \xi_0=x\Big) \, d{\tilde\mu}_{k,j}(x)\\
   &<\sum_{k,j} p_k p_j \int\limits_0^\infty x\IE_{k,j} \Big(
\frac{\C(\T_0^*)}{\C(\T)} 
\mid \xi_0=\infty\Big) \, d{\tilde\mu}_{k,j}(x)\\
    &=\sum_{k,j} p_k p_j  E\Big(\!\frac{Z_1^{(j)}+\cdots
+Z_j^{(j)}}{Z_1^{(j)}+\cdots +Z_j^{(j)}
 +Z_1^{(k)}+\cdots +Z_k^{(k)}}\!\Big) \!\int\limits_0^\infty x \, 
d{\tilde\mu}_{k,j}(x)\\
&=\gamma \sum_{k,j}r_{k,j} E\Big(\frac{Z_1^{(j)}+\cdots
+Z_j^{(j)}}{Z_1^{(j)}+\cdots +Z_j^{(j)}
 +Z_1^{(k)}+\cdots +Z_k^{(k)}}\Big).
\end{align*}
To see that the inequality in the above calculation is strict,
observe that $\C(\T)=\C(\T_0^*)+\cdots+\C(\T_k^*)$ on
$\{\ii(\0)=k\}$, and when the conductance of $w_0(\0)$ increases,
so does the effective conductance of~$\T_0^*$, and therefore so
does the quantity $\frac{\C(\T_0^*)}{\C(\T)}$; note also that
putting an infinite conductance to an edge means effectively
shrinking this edge.
Hence, due to \eqref{sum_pkpj},
\begin{equation}
\label{domhalf}
\IE\Big(\xi_0 \frac{\C(\T_0^*)}{\C(\T)} \Big) < \frac{\gamma}{2}\, .
\end{equation}
Thus, with~\eqref{eq:speed'} and~\eqref{domhalf}, we obtain
$v > 0$,
which concludes the proof of Theorem~\ref{pos_speed}. 
\qed
\noindent
\textit{Proof of Theorem~\ref{compara_speeds}.} Observe that, by 
symmetry,
\begin{equation}
\label{annealed_escape}
 \IE_k \Big(\frac{\C(\T_0^*)}{\C(\T)}\Big) = \frac{1}{k+1}.
\end{equation}
So, from~(\ref{speed_SRW}) we obtain that
\[
 v_{SRW} = 1 - \frac{2}{\gamma} \IE\left(\xi_0\right)\,
\IE\Big(\frac{\C(\T_0^*)}{\C(\T)}\Big),
\]
and~\eqref{speed_cov} follows from~\eqref{eq:speed'}.
Let us now prove part~(ii).
 From~\eqref{annealed_escape}
we obtain that
\[
 \sum_m p_m \IE_{k,m} \Big(\frac{\C(\T_0^*)}{\C(\T)}\Big) = 
\frac{1}{k+1}.
\]
Once again, we observe that, when $\xi_0$
increases (while fixing the other conductances), so does
$\frac{\C(\T_0^*)}{\C(\T)}$;
 this means that $\xi_0$ and $\C(\T_0^*)/\C(\T)$ are
positively correlated under $\IE_{k,m}$ (and strictly positively
correlated for at least one pair $(k,m)$ in the case when~$\xi_0$ is
a nondegenerate random variable), so we have
\begin{align*}
  \IE \Big(\xi_0 \frac{\C(\T_0^*)}{\C(\T)}\Big) &=
  \sum_{k,m} p_kp_m\IE_{k,m}\Big(\xi_0
\frac{\C(\T_0^*)}{\C(\T)}\Big)\\
 & > \sum_{k,m} p_kp_m \gamma_{k,m}
\IE_{k,m}\Big(\frac{\C(\T_0^*)}{\C(\T)}\Big)\\
 & =\gamma \sum_k \frac{p_k}{k+1}\\
 & = \IE\left(\xi_0\right) \, 
\IE\Big(\frac{\C(\T_0^*)}{\C(\T)}\Big),
\end{align*}
where we used $\gamma_{k,m} = \gamma$, for all $k,m$ 
for the third equality. Now part~(ii) follows
from~\eqref{speed_cov}.\qed

\noindent 
\textit{Proof of Theorem~\ref{zero_speed}}.
Assume that there exist $k, m$ such that $\gamma_{k,m} = \infty$. 
We will show that for $T_n: = \inf\{j: |X_j| = n\}$,
\begin{equation}
\label{Tgrows}
\frac{T_n}{n} \to \infty \quad \text{$\Pomega$-a.s. for $\IP$-almost
all~$\omega$.}
\end{equation}
Since, for any $\varepsilon > 0$, 
$\{|X_n| \geq  \lfloor n\varepsilon\rfloor\} 
\subseteq \{T_{\lfloor n\varepsilon\rfloor} \leq n\}$, 
\eqref{Tgrows} implies that 
$\frac{|X_n|}{n} \to 0$, $\Pomega$-a.s. for $\IP$-almost all
$\omega$. 
To show~\eqref{Tgrows}, we will prove that there is an 
i.i.d.\ sequence of
random variables 
$(\eta_j)_{j \geq 1}$ with infinite expectations such
that $T_n$ is larger than $\frac{1}{\lfloor
n/5 \rfloor}\sum_{i=1}^{\lfloor n/5 \rfloor}\eta_i$.
Roughly speaking, the infinite expectations come from the fact that
the random walk frequently crosses
bonds~$(y,z)$
with $\ii(y)=k$ and $\ii(z)=m$, where the conductances of the
neighbouring bonds are not too large.
To understand the following proof, it is good to keep in mind 
that we can 
construct the tree successively with the random walk, adding new 
vertices and edges as the random walk explores the tree.

Let $M,C>0$ (to be specified later). For any $x\neq \0$ we denote 
by $\overleftarrow x$ the predecessor vertex with respect to~$x$,
i.e., $\overleftarrow x$ is the neighbor of~$x$ such that 
$|\overleftarrow x|=|x|-1$.
Let a vertex~$x\neq \0$ be \textit{good} if $\ii(x)\leq M$,
$\ii(\overleftarrow x) \leq M$, and $\xi(\overleftarrow
x,x)\leq C$, i.e.,\ the bond from~$x$
towards the root has conductance at most~$C$, while the degrees
of~$x$ and its predecessor are not too large. 

\begin{figure}
\centering
\includegraphics{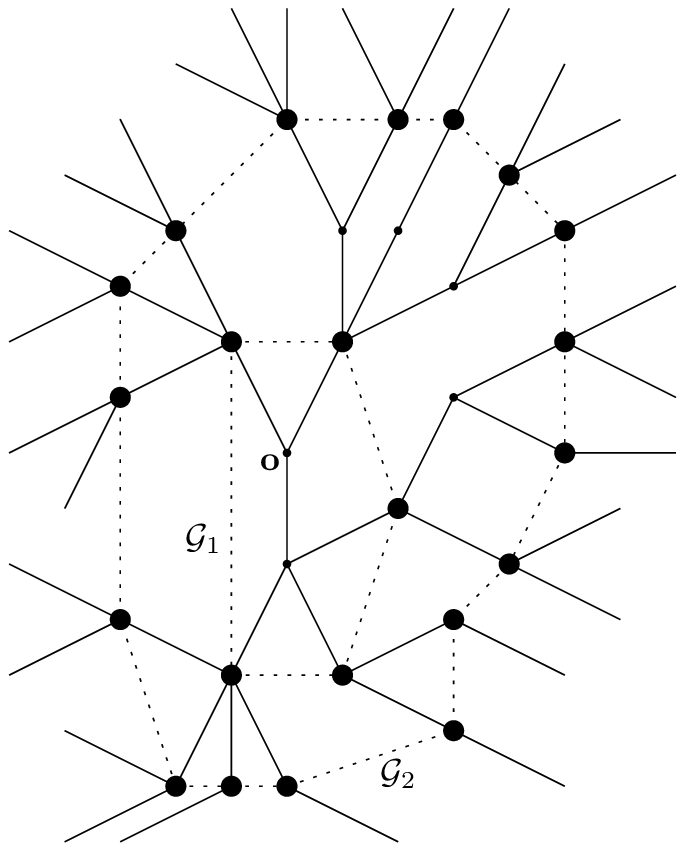}
\caption{On the definition of the cutsets $\G_1,\G_2,\ldots$ (good
sites are marked by larger circles)} 
\label{f_cutsets}
\end{figure}

We now define recursively cutsets of good vertices which the random
 walk has to cross on its way. 
For $u, v \in \T$ with $u<v$,  let a \lq\lq ray from~$u$ to
$v$\rq\rq\
be a path $(z_1, \ldots , z_K)$, with $z_1 = u$ and 
$|z_{i+1}| = |z_i| +1$, $\forall i$ and $z_K = v$. 
Call a vertex \textit{bad} if it is not good.
Let $\G_1$ be the set of all vertices $u_1$ which are good and such 
that all vertices on the ray from the root to $u_1$ 
are bad.
 Then, let $\G_2$ be the set of all vertices $u_2$ which are good
and
 such that the ray from the root to $u_2$ contains exactly one good
vertex $u_1 \in \G_1$ with $|u_1| < |u_2|$, and so on, see
Figure~\ref{f_cutsets}. 
Let $B_n : =
\{u \in \T: |u| \leq n\}$.

\medskip
\noindent
\textbf{Claim.} We can choose large enough $M,C$ in 
such a way that
\begin{equation}
\label{eq_claim}
 \IP[\G_{\lfloor n/5\rfloor } \subseteq B_n \hbox { for all } n 
\hbox{ large enough}] =1\, .
\end{equation}

\medskip
\noindent
\textit{Proof of the claim.}
If $\G_{\lfloor n/5\rfloor } \not\subseteq B_n$, there has to
be a ray from the root to a vertex at distance~$n$ from the root,
containing at least $4n/5$ bad vertices; we will show that this
happens with exponentially small probability and so one 
obtains~\eqref{eq_claim} from the Borel--Cantelli lemma. 

First, let us prove that, for large enough~$M$, with large
probability on every path to the level~$n$ there are at most~$n/5$
sites with index greater than~$M$. For this, consider a branching
random walk starting with one particle at the origin, described in
the following way:
\begin{itemize}
 \item on the first step the particle generates $j+1$ offspring with
probabilities $p_j$, $j\geq 1$, and on subsequent time moments every
particle generates~$j$ offspring with probabilities $p_j$, $j\geq
1$, independently of the others;
 \item if the number of a particle's offspring is less than
or equal to~$M$, then all
the offspring stay on the same place, and if it is greater than~$M$,
then all its offspring go one unit to the right.
\end{itemize}
With this interpretation, we have to
prove that with large probability at time~$n$ the whole cloud is
to the left of~$n/5$. 
In fact, it is well-known that the position of the
righthmost particle grows linearly in time, and the linear speed goes
to~$0$ if~$M$ goes to~$\infty$; there are several possible ways to
show this.  For instance, one can
use the many-to-one lemma (see e.g.\ formula~(2.2)
of~\cite{AidShi:10}),
dealing with the small difficulty that at time~$1$ the offspring
distribution is different. Another possibility is to consider
the process 
\[
 Z_n = (2\mu)^{-n}\sum_{k\in\Z}\eta_n(k)(2\mu)^{6k},
\]
where $\mu = \sum_{j=1}^\infty j p_j$ and $\eta_n(k)$ is the number
of particles of the branching
random walk at time~$n$ at site~$k$. With a straightforward
calculation, one obtains that if~$M$ is large enough, then~$Z$
is a (nonnegative) supermartingale. So, we obtain
\begin{align}
 \IP[\text{there exists }k\geq n/5 \text{ such that }\eta_n(k)\geq 1]
  & \leq \IP[Z_n\geq (2\mu)^{-n}\cdot (2\mu)^{6n/5}]\nonumber\\
 &\leq \frac{\IE Z_n}{(2\mu)^{n/5}}, \nonumber\\
 &\leq (2\mu)^{-n/5}\, , \label{cloud_particles}
\end{align}
using in the last inequality  the fact
that~$Z$ is a supermartingale.

Now, if every path to the level~$n$ contains at most~$n/5$ sites
with index greater than~$M$, then on every path to the level~$n$
there are at least $3n/5-1$
sites with index less than or equal to~$M$ and such the predecessor
site has index less than or equal to~$M$ as well. Also, using
the Chebychev inequality one immediately obtains that
with probability at least $1-2^{-n}$
the total number of paths to level~$n$ is less than~$(2\mu)^n$. 
Next, denoting by
\[
 h(C) = \max_{i,j\leq M} {\tilde\mu}_{i,j}(C, +\infty),
\]
we have, clearly, that $h(C)\to 0$ as $C\to \infty$.
Let us choose~$C$
in such a way that~$h(C)$ is small enough to assure the following: 
on a fixed path to level~$n$ (with given degrees of vertices but the
conductances not yet chosen) such that the number of bonds there
that belong to $\cup_{i,j\leq M}U_{i,j}$ is at least $3n/5-1$, 
 the number of good sites is at least~$n/5$ with probability 
at least $1-(3\mu)^{-n}$ (this amounts to estimating the
probability that a sum of $3n/5-1$ Bernoulli random variables with
probability of success $1-h(C)$ is at least $n/5$). Then, we use the
union bound and the Borel--Cantelli lemma to conclude the proof of the
claim.
\qed

Now, define by ${\tilde T}_j=\min\{n:X_n\in\G_{4j}\}$,
$k=1,2,3,\ldots$, the hitting times of the sets
$\G_4,\G_8,\G_{12},\ldots$ (for formal reasons, we also set ${\tilde
T}_0:=0$). Without restricting generality, one can assume that
$M\geq\max\{k,m\}$ (recall that $k,m$ are such that
$\gamma_{k,m}=\infty$). 
Consider the events $A_j$, $j=1,2,3,\ldots$, defined in
the following way:
\begin{align*}
 A_j &= \Big\{\text{there exist }y,z \text{ with }\ii(y) = k,
\ii(z)= m,
\text{ such that }
 y=\overleftarrow z, X_{{\tilde T}_j}=\overleftarrow y, \\
& ~~~~~~~
\text{ and } C^{-1} \leq \xi(e) \leq C \text{ for all }
e\neq (y,z) \text{ such that } e\sim y \text{ or }e\sim z\Big\}
\end{align*}
(observe that the event~$A_j$ concerns the yet unexplored part of
the tree at time ${\tilde T}_j$).
Note that there is some $g = g(M,C) > 0$ such that
\begin{equation}
\label{gotoy}
\IE\Po[A_j] \geq g(M,C)\, .
\end{equation}
Further, if $X_{{\tilde T}_j} = x$, then, since $x$ is good, the
probability to go from~$x$ to~$y$ (i.e., the site in the definition
of the event $A_j$) is bounded below by $\frac{1}{(1+M)C^2}$.
Then, the number of subsequent crossings $N_{(y,z)}$ 
of the bond 
$(y,z)$ (again, $y,z$ are the sites from the definition
of the event $A_j$) dominates a geometric
random variable with parameter $h_0:=\frac{(k+m)C}{\xi(y,z)+(k+m)C}$.
So, under the averaged measure $\IE\Po=\int \Po[\,\cdot\,]
\IP(d\omega)$, each of the random variables
$({\tilde T}_j-{\tilde T}_{j-1})$ dominates a random variable
$\eta_j$ with law
\[
 \eta_j = \begin{cases}
           0, & \text{with probability } 1-\frac{g(M,C)}{(1+M)C^2},\\
           \mathsf{Geometric}(h_0), & \text{with probability }
\frac{g(M,C)}{(1+M)C^2},
          \end{cases}
\]
and $\eta_1,\eta_2,\eta_3,\ldots$ are i.i.d.\ under the
measure $\IE\Po$. Since, clearly, the expectation of~$\eta_1$ under
the averaged measure $\IE\Po$ is infinite,
this implies
Theorem~\ref{zero_speed} as explained in the beginning of the proof.
\qed

\begin{figure}
\centering
\includegraphics{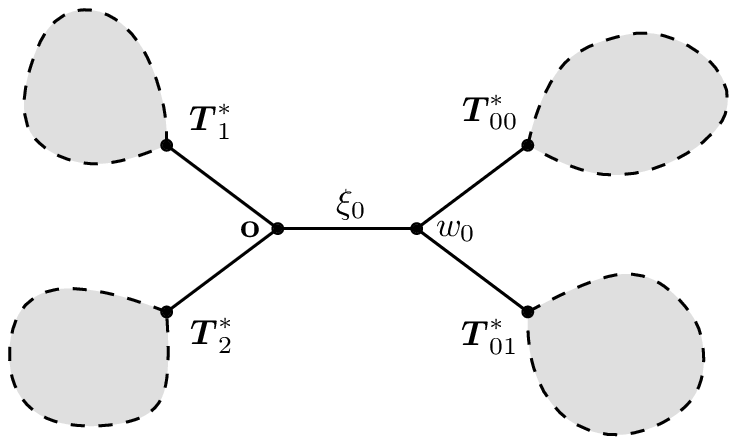}
\caption{On the definition of the trees~$\T_{00}^*, \T_{01}^*$}
\label{f_T00}
\end{figure}

\noindent \textit{Proof of Proposition~\ref{ex1}}. For the binary
tree, equation~\eqref{eq:speed'} implies that (see
Figure~\ref{f_T00})
\begin{align}
 v&=1-\frac{2}{\gamma}\IE \Big(\xi_0
\frac{\C(\T_0^*)}{\C(\T)}\Big)
\label{gg1}\\
&=1-\frac{2}{\gamma}\IE\Big(
 \xi_0\frac{\big(1+\frac{\C(\T_{00}^*)+
\C(\T_{01}^*)}{\xi_0}\big)^{-1}(\C(\T_{00}^*)+\C(\T_{01}^*))}
{\C(\T_0^*)+\C(\T_1^*)+\C(\T_2^*)}\Big).
\label{gg2}
\end{align}
Then, we can write
\begin{equation}
\label{aeps1}
 \IE \Big(\xi_0\frac{\C(\T_0^*)}{\C(\T)}\Big)=
  (1-\eps_n)\IE \Big(\frac{\C(\T_0^*)}{\C(\T)}\mid \xi_0=1\Big)
+ \eps_n a_n\IE \Big(\frac{\C(\T_0^*)}{\C(\T)}\mid \xi_0=a_n\Big).
\end{equation}
Also, by symmetry we have
\begin{align}
 \frac{1}{3}&=\IE \left(\frac{\C(\T_0^*)}{\C(\T)}\right)\nonumber\\
&=(1-\eps_n)\IE \Big(\frac{\C(\T_0^*)}{\C(\T)}\mid \xi_0=1\Big) +
\eps_n\IE \Big(\frac{\C(\T_0^*)}{\C(\T)}\mid \xi_0=a_n\Big).
\label{aeps2}
\end{align}
Since~$\C(\T_0^*)/\C(\T) \le 1$, we obtain from~\eqref{aeps2} that
\begin{equation}
\label{ll1}
 \IE \Big(\frac{\C(\T_0^*)}{\C(\T)}\mid \xi_0=1\Big) 
        \to \frac{1}{3} \qquad \text{as $n\to\infty$}
\end{equation}
(observe that the expectation in the left-hand side depends on~$n$
in fact), and so, by~\eqref{gg1} and~\eqref{aeps1}, we have
$v_n\to 1-2/3=1/3$ in the case $a_n \eps_n \to 0$ (note that in
this case~$\gamma =1-\eps_n+a_n \eps_n\to 1$).

Now we consider the two other cases. First, we want to show that
\[
\IE \Big(\frac{\C(\T_0^*)}{\C(\T)}\mid \xi_0=a_n\Big)\to
\frac{1}{2} \qquad \text{as $n\to\infty$}.
\]
Putting an infinite conductance  to the edge $\ell_0$, we
obtain (as Figure~\ref{f_T00} suggests)
\[
\C(\T)< \C(\T_{00}^*)+\C(\T_{01}^*)+\C(\T_1^*)+\C(\T_2^*)
\]
(naturally, $w_0$ is supposed to be the root of $\T_{00}^*$ 
and $\T_{01}^*$). Then,
\begin{align*}
 \lefteqn{\IE \Big(\frac{\C(\T_0^*)}{\C(\T)}\mid \xi_0=a_n\Big)}&\\
   &>
    \IE \Big(\frac{\C(\T_0^*)}{ \C(\T_{00}^*)+\C(\T_{01}^*)+
\C(\T_1^*)+\C(\T_2^*) } \mid \xi_0=a_n\Big)\\
&=\frac{1}{2}
  -\IE \Big(\frac{\C(\T_{00}^*)+\C(\T_{01}^*)
    -\big(\frac{1}{a_n}+\frac{1}{\C(\T_{00}^*)
       +\C(\T_{01}^*) }\big)^{-1}}{
       \C(\T_{00}^*)+\C(\T_{01}^*)+\C(\T_1^*)+\C(\T_2^*) }\mid 
      \xi_0=a_n\Big)\\
&=\frac{1}{2}
  -\IE \Big(\frac{\C(\T_{00}^*)+\C(\T_{01}^*)}
{ \C(\T_{00}^*)+\C(\T_{01}^*)+\C(\T_1^*)+\C(\T_2^*)}\\
&\qquad\qquad\qquad \times\Big[1-\Big(\frac{
\C(\T_{00}^*)+\C(\T_{01}^*)
}{a_n}+1\Big)^{-1} \Big] \Big).
\end{align*}
Observe that
\[
\frac{\C(\T_{00}^*)+\C(\T_{01}^*)}{
\C(\T_{00}^*)+\C(\T_{01}^*)+\C(\T_1^*)+\C(\T_2^*)}< 1,
\]
and (because if the conductance on the first edge is~$1$, then the
effective conductance of the tree is less than~$1$)
\[
\IP[ \C(\T_{00}^*)\le 1]\ge 1-\eps_n.
\]
Thus, we have that
\[
 \frac{ \C(\T_{00}^*)+\C(\T_{01}^*) }{a_n}\to 0
\]
in probability and so
\[
 \frac{\C(\T_{00}^*)+\C(\T_{01}^*)}{ \C(\T_{00}^*)
  +\C(\T_{01}^*)+\C(\T_1^*)+\C(\T_2^*)}
\Big( 1-\Big(\frac{ \C(\T_{00}^*)+\C(\T_{01}^*) }{a_n}+1\Big)^{-1}
\Big) \to 0
\]
in probability and hence in~$L_1$. Thus, we indeed have
\begin{equation}
\label{ll2}
 \IE \Big(\frac{\C(\T_0^*)}{\C(\T)}\mid \xi_0=a_n\Big)\to 
        \frac{1}{2}.
\end{equation}

When~$a_n\eps_n\to \infty$, using~\eqref{aeps1}, \eqref{ll1},
\eqref{ll2},
 we obtain
\[
\frac{1}{\gamma}\IE \Big(\xi_0\frac{\C(\T_0^*)}{\C(\T)}\Big)\to
\frac{1}{2}
\]
and so~$v_n\to 0$ by~\eqref{gg1}.

When~$a_n\eps_n\to\eta\in(0,\infty)$, we have by~\eqref{aeps1},
\eqref{ll1}, \eqref{ll2}, that
\[
\frac{1}{\gamma}\IE \Big(\xi_0\frac{\C(\T_0^*)}{\C(\T)}\Big)
    \to \frac{1}{1+\eta}\Big(\frac{1}{3}+\eta\frac{1}{2}\Big)
\]
as $n\to\infty$, and so
\[
 v_n\to 1-\frac{2}{1+\eta}\Big(\frac{1}{3}+
\frac{\eta}{2}\Big) =\frac{1}{3(\eta+1)},
\]
which finishes the proof of Proposition~\ref{ex1}. \qed

\section{Open questions}
\label{s_open}

\begin{itemize}
\item[1.]
We conjecture that (ii) in Theorem~\ref{compara_speeds} 
still holds in the case where $\gamma < \infty$ and the
$\gamma_{k,m}$'s are different. This amounts to proving that
\[
 \Cov\Big(\xi_0,\frac{\C(\T^*_0)}{\C(\T)}\Big) \geq 0
\]
for this case.
\item[2.] 
If $\gamma_{k,m} < \infty$ for all $k,m$ but $\gamma = \infty$, it is
not clear under which conditions the speed of the random walk is zero
or strictly positive, respectively. We believe that both can happen.
\item[3.]

\textbf{Problem:} Find  conditions  for graphs on which the SRW has positive speed
such that  for the random conductance model, taking i.i.d.\ conductances with
finite mean, the speed of the corresponding random walk is less or
equal, or strictly less than the speed of SRW.

\item[4.] As mentioned in the introduction,  our random conductance model can be seen as a 
unimodular random network, under the condition that $\gamma<\infty$. This suggests to formulate an interesting special case of the
above problem:

\textbf{Question:} Is it true that all non-amenable unimodular random
graphs exhibit the \emph{slowdown} phenomenon, i.e. that for the random
conductance model, taking (non-degenerate)  i.i.d.\ conductances with finite
mean, the speed of the corresponding random walk is strictly less than the
speed of the simple random walk?

\end{itemize}

\section*{Acknowledgements}
S.M.\ was partially supported by FAPESP (2009/08665--6).
S.P.\ was partially supported by  CNPq (300328/2005--2). 
M.V.\ was partially supported 
         by CNPq (304561/2006--1). S.P. and M.V.
thank FAPESP (2009/52379--8) for financial support.
The work of N.G.\ was partially supported by FAPESP (2010/16085--7).
We also thank CAPES/DAAD (Probral)
for support.


\end{document}